\documentclass[12pt,letter]{article}
\usepackage{amsfonts}

\usepackage{graphicx}
\usepackage{amsmath}
\usepackage{citesort}

\begin{document}

\begin{center}
{\textbf{\Large A series solution and a fast algorithm\\
for the inversion of \\[2mm]
the spherical mean Radon transform}}\\[6mm]

L. Kunyansky \\[2mm]

University of Arizona, Tucson
\end{center}

\medskip

\begin{abstract}
An explicit series solution is proposed for the inversion of the spherical
mean Radon transform. Such an inversion is required in problems of thermo-
and photo- acoustic tomography. Closed-form inversion formulae are currently
known only for the case when the centers of the integration spheres lie on a
sphere surrounding the support of the unknown function, or on certain unbounded
surfaces. Our approach results
in an explicit series solution for any closed measuring surface surrounding
a region for which the eigenfunctions of the Dirichlet Laplacian are
explicitly known --- such as, for example, cube, finite cylinder,
half-sphere etc. In addition, we present a fast reconstruction algorithm
applicable in the case when the detectors (the centers of the integration
spheres) lie on a surface of a cube. This algorithm reconsrtucts 3-D images
thousands times faster than backprojection-type methods.
\end{abstract}

\section*{Introduction}

The problem of image reconstruction in thermo-acoustic and photo-acoustic
tomography is equivalent to recovering a function from a certain set of its
spherical means \cite{kruger1,kruger,XuWang,XuWang0,kuchrev}.
The process starts with object of interest being excited by
a short electromagnetic pulse. This causes
thermal expansion of the tissue, and generates an acoustic wave whose
intensity is recorded by a set of detectors located outside the object.
The intensity of the thermal expansion depends on the local properties of
the tissue, and is of interest to a doctor, since abnormally high values of
this function are indicative of a tumor. Under certain simplifying
assumptions, the measurements can be related to the integrals of the
expansion intensity over the spheres with the centers at the detectors'
locations. The reconstruction of the local properties from these integrals
is equivalent to the inversion of the spherical mean Radon transform.

Some of the recent results on the injectivity of this transform as well as
the corresponding range conditions can be found in \cite
{AQ,trinity,Finch,AK1,Finch1,AK2}. In the present paper we concentrate on
inversion formulae and algorithms for the solution of the reconstruction
problem. Generally, in such applications as photo- and thermo- acoustic
tomography, the designer of the measuring system has a certain freedom
of choice when selecting
the detectors' positions (the centers of the integration spheres). Most
of the known explicit solutions pertain to the spherical acquisition
geometry, in other words to the configuration in which the detectors are
located on a sphere surrounding the object. Such are the recently found
series solutions~\cite{Norton2D,Nort3D,XuWang0} and backprojection-type
formulae~\cite{Finch,XuWang,Finch2,kunya}. Explicit reconstruction
formulae are also known for such acquisition geometires as an infinite
plane~\cite{Natt,XuWang,Faw} and an infinite cylinder~\cite{XuWang}.

The spherical geometry is preferable to unbounded measuring surfaces since
the latter have to be truncated in practice, which leads to errors in the
reconstruction. However, there are compelling reasons to consider other
non-spherical bounded measuring surfaces as well. For example, as shown in
section~2, if the detectors are located on a surface of a cube surrounding
the object of investigation, it is possible to design a fast algorithm that
reconstructs the unknown 3-D function in a matter of seconds --- as opposed
to several hours required for the algorithms based on straightforward
discretization of one of the 3-D backprojection-type inversion
formulae~\cite{Finch,XuWang,kunya}.


We thus present a series solution for the inversion of spherical mean Radon
transform in the case when the centers of the integration spheres lie on a
closed surface surrounding a bounded connected region in
$\mathbb{R}^{n},\quad n\geq 2$.
Our procedure requires knowledge of the eigenfunctions
of the Dirichlet Laplacian defined on the region enclosed by the
measuring surface. For many regions of practical interest such
eigenfunctions are known explicitly. Among such regions in 3-D are, for
example, a rectangle, a ball, an ellipsoid, a cylinder, a spherical shell.
In addition, these eigenfunctions can be easily found for certain subsets of
these bodies obtained by dissecting them along a plane of symmetry --- for
example for a half-ball, half-cylinder, certain triangular prisms and
tetrahedra. Yet another example of regions with explicitly known
eigenfunctions is given by the crystallographic domains (see \cite
{berard1,berard2} for details). A generalization of this approach
to a general connected region is possible if one computes
the eigenfunctions of the Dirichlet Laplacian  numerically.
In this case, however, the reconstruction algorithm is likely to
be rather expensive from the computational point of view.

The proof of the range theorem in~\cite{trinity} involves implicitly a
reconstruction procedure also based on eigenfunction expansions. Unlike the
present method, that procedure would involve division of analytic functions
that have countable number of zeros. While the range theorem guarantees
cancellation of these zeros when the data are in the range of the direct
transform, a stable numerical implementation of such division would be
complicated if not impossible. (A similar problem arises with the series
solution of~\cite{Norton2D} that involves division of certain computed
quantities by Bessel functions.) The technique we present below does not
require such divisions.

Section~1 contains a general description of the present method.
The efficiency of numerical realization of this technique depends,
in particular, on the availability of fast algorithms for the
summation of the arising eigenfunction expansions.
In the simplest case of a cubic (or rectangular) measuring surface such an
algorithm is the 3-D Fast Sine Fourier transform. This allows us to design a
very efficient reconstruction algorithm for this particular measuring
configuration, as discussed in section~2. Finally, in section~3 we
investigate an interesting property that seems to be exclusive to the
series solutions presented in this paper. Namely, this technique
 will produce a theoretically exact image within the region enclosed
by the measuring surface even if there are sources outside that region.
This property can prove to be useful for reducing the sensitivity
of the measuring system to external noise. Such a noise cancellation
will occur, however, only if all the measurements are performed
simultaneously by a fixed set of detectors; a synthesized measuring
surface will not exhibit this phenomenon.

\section{Series solution\label{formulation}}

Suppose that $C_{0}^{1}$ function $f(\mathbf{x})$,
$\mathbf{x\in }\mathbb{R}^{n},\quad n\geq 2$
is compactly supported within the bounded connected open
region $\Omega $ with boundary $\partial \Omega .$ Our goal is to
reconstruct $f(\mathbf{x})$ from its projections $g(\mathbf{z},r)$ defined
as the integrals of $f(\mathbf{x})$ over the spheres of radius $r$ centered
at $\mathbf{z}$:

\begin{equation*}
g(\mathbf{z},r)=\int\limits_{\mathbb{S}^{n-1}}f(\mathbf{z}+r\hat{t})r^{n-1}
ds(\hat{t}\mathbf{),}
\end{equation*}
where $\mathbb{S}^{n-1}$ is the unit sphere in $\mathbb{R}^{n},$ $\hat{t}$
is a unit vector, and $ds$ is the normalized measure in $\mathbb{R}^{n}$.
Projections are assumed to be known for all $\mathbf{z}\in \partial \Omega ,$
$0\leq r\leq \mathrm{diam(}\Omega )$ (integrals for $r>\mathrm{diam(}\Omega )
$ automatically equal zero, since the corresponding integration spheres do
not intersect the support of the function).

Suppose $\lambda _{m}^{2}$, $u_{m}(\mathbf{x})$ are the eigenvalues and
normalized eigenfunctions of the Dirichlet Laplacian $- \Delta$ on $\Omega $
with zero boundary conditions, i.e.
\begin{align}
\Delta u_{m}(\mathbf{x})+\lambda _{m}^{2}u_{m}(\mathbf{x})& =0,\qquad
\mathbf{x}\in \Omega ,\quad \Omega \subseteq \mathbb{R}^{n},  \label{Helmeq}
\\
u_{m}(\mathbf{x})& =0,\qquad \mathbf{x}\in \partial \Omega ,  \notag \\
||u_{m}||_{2}^{2}& \equiv \int\limits_{\Omega }|u_{m}(\mathbf{x})|^{2}
d\mathbf{x}=1.  \notag
\end{align}
We would like to reconstruct function $f(\mathbf{x})$ from the known values
of its spherical integrals $g(\mathbf{z},r)$ with the centers on $\partial
\Omega $:
\begin{equation*}
g(\mathbf{z},r)=\int\limits_{\mathbb{S}^{n-1}}f(\mathbf{z}+r\mathbf{\hat{s}}
)r^{n-1}d\mathbf{\hat{s}},\qquad \mathbf{z}\in \partial \Omega .
\end{equation*}
We notice that $u_{m}(\mathbf{x})$ is the solution of the Dirichlet problem
for the Helmholtz equation with zero boundary conditions and the wave number
$\lambda _{m}$, and thus it admits the Helmholtz representation
\begin{equation}
u_{m}(\mathbf{x})=\int_{\partial \Omega }\Phi _{\lambda _{m}}(|\mathbf{x}-
\mathbf{z|})\frac{\partial }{\partial \mathbf{n}}u_{m}(\mathbf{z})ds(
\mathbf{z)}\qquad \mathbf{x}\in \Omega ,  \label{helmdiscr}
\end{equation}
where $\Phi _{\lambda _{m}}(|\mathbf{x}-\mathbf{z|})$ is a free-space
rotationally invariant Green's function of the Helmholtz equation~(\ref
{Helmeq}).

Our approach is based on the fact that eigenfunctions $\left\{ u_{m}
(\mathbf{x})\right\} _{0}^{\infty }$
form an orthonormal basis in $L_{2}(\Omega ).$
Therefore $f(\mathbf{x})$ can be represented by the series
\begin{equation}
f(\mathbf{x})=\sum_{m=0}^{\infty }\alpha _{m}u_{m}(\mathbf{x})
\label{fourierser}
\end{equation}
with
\begin{equation}
\alpha _{m}=\int_{\Omega }u_{m}(\mathbf{x})f(\mathbf{x})d\mathbf{x.}
\label{serkoef}
\end{equation}
Since $f(\mathbf{x})$ is $C_{0}^{1},$ series~(\ref{fourierser}) converges
pointwise. The reconstruction formula will result if we substitute
representation~(\ref{helmdiscr}) into (\ref{serkoef}) and interchange the
order of integrations
\begin{align}
\alpha _{m}& =\int_{\Omega }u_{m}(\mathbf{x})f(\mathbf{x})d\mathbf{x}  \notag
\\
& =\int_{\Omega }\left( \int_{\partial \Omega }\Phi _{\lambda _{m}}
(|\mathbf{x}-\mathbf{z}|)
\frac{\partial }{\partial \mathbf{n}}u_{m}(\mathbf{z})ds
(\mathbf{z)}\right) f(\mathbf{x})d\mathbf{x}  \notag \\
& =\int_{\partial \Omega }\left( \int_{\Omega }\Phi _{\lambda _{m}}
(|\mathbf{x}-\mathbf{z}|)f(\mathbf{x})d\mathbf{x}\right)
\frac{\partial }{\partial\mathbf{n}}
u_{m}(\mathbf{z})ds(\mathbf{z)}  \label{serkoef1} \\
& =\int_{\partial \Omega }I(\mathbf{z},\lambda _{m})\frac{\partial }
{\partial \mathbf{n}}u_{m}(\mathbf{z})ds(\mathbf{z),}  \label{serkoef1a}
\end{align}
where
\begin{equation*}
I(\mathbf{z},\lambda _{m})=\int_{\Omega }\Phi _{\lambda _{m}}(|\mathbf{x}-
\mathbf{z}|)f(\mathbf{x})d\mathbf{x.}
\end{equation*}
The change of the integration order is justified by the continuity of
eigenfunctions $u_{m}(\mathbf{x})$. Function $I(\mathbf{z},\lambda _{m})$ is
easily computed from the projections
\begin{equation*}
I(\mathbf{z},\lambda _{m})=\int_{\Omega }\Phi _{\lambda _{m}}(|\mathbf{x}-
\mathbf{z}|)f(\mathbf{x})d\mathbf{x}=\int\limits_{\mathbb{R}^{+}}
g(\mathbf{z},r)\Phi _{\lambda _{m}}(r)dr,
\end{equation*}
and with Fourier coefficients $\alpha _{m}$ now known, $f(\mathbf{x})$ is
reconstructed by summing series~(\ref{fourierser}).

If desired, this solution can be re-written in the form of a
backprojection-type formula:
\begin{align}
f(\mathbf{x})& =\sum_{m=0}^{\infty }\alpha _{m}u_{m}(\mathbf{x})
=\int_{\partial \Omega }\left( \sum_{m=0}^{\infty }\alpha _{m}\Phi
_{\lambda _{m}}(|\mathbf{x}-\mathbf{z}|)\frac{\partial }{\partial \mathbf{n}}
u_{m}(\mathbf{z})\right) ds(\mathbf{z)}  \notag \\
& =\int_{\partial \Omega }h(\mathbf{z},|\mathbf{x}-\mathbf{z|})ds(\mathbf{z),}
\label{serbackpr}
\end{align}
where
\begin{equation}
h(\mathbf{z},t)=\sum_{m=0}^{\infty }\alpha _{m}\Phi _{\lambda _{m}}(t)
\frac{\partial }{\partial \mathbf{n}}u_{m}(\mathbf{z}),  \label{serfiltr}
\end{equation}
and coefficients $a_{m}$ are computed using equation~(\ref{serkoef1a}). In
the above formula equation~(\ref{serbackpr}) is clearly a backprojection
operator, and~(\ref{serfiltr}) is a filtration. However, the latter operator
is now represented by a series rather than by a closed form expression.
Moreover, this operator is not local in $\mathbf{z,}$ unlike the filtration
operator of the known closed-form explicit inversion
formulae~\cite{Finch,Finch2,XuWang,kunya}.

Finally, we notice that if function $f(\mathbf{x})$ is not smooth but rather
belongs to $L^{2}(\Omega ),$ our reconstruction formulae are still valid if
equation~(\ref{fourierser}) is understood in the $L^{2}$ sense.

\section{A fast algorithm for the cubic measurement surface in 3D}

A cube is the simplest of the regions whose eigenfunctions of the Dirichlet
Laplacian are known explicitly; they are products of sine functions. In the
present section we exploit the simple structure of these eigenfunctions to
develop a fast reconstruction algorithm applicable in the case when the
detectors are located on a surface of a cube (a generalization to a
rectangular case is straightforward). Such a measuring surface can be either
sampled by regular detectors or synthesized from measurements made by
interferometric line detectors as discussed in the Introduction.

Let the sought function $f(\mathbf{x})$ be supported within the cube $\Omega
=[0,R]\times [0,R]\times [0,R].$ We will index the normalized
eigenfunctions $u_{\mathbf{m}}(\mathbf{x})$ and eigenvalues
$\lambda _{\mathbf{m}}$
of the Dirichlet Laplacian on this region using vector
$\mathbf{m}=(m_{1},m_{2},m_{3}),$ $m_{1},m_{2},m_{3}\in \mathbb{N:}$
\begin{eqnarray*}
u_{\mathbf{m}}(\mathbf{x}) &=&\frac{8}{R^{3}}\sin \frac{\pi m_{1}x_{1}}{R}
\sin \frac{\pi m_{2}x_{2}}{R}\sin \frac{\pi m_{3}x_{3}}{R}, \\
\lambda _{\mathbf{m}} &=&\pi ^{2}|\mathbf{m|}^{2}.
\end{eqnarray*}
Cube $\Omega $ has six faces $\delta \Omega _{i},i=1,...,6:$
\begin{eqnarray*}
\delta \Omega _{1} &=&\{\mathbf{x}|x_{1}=R,0<x_{2}<R,0<x_{3}<R,\}, \\
\delta \Omega _{2} &=&\{\mathbf{x}|x_{1}=0,0<x_{2}<R,0<x_{3}<R,\}, \\
\delta \Omega _{3} &=&\{\mathbf{x}|x_{2}=R,0<x_{1}<R,0<x_{3}<R,\}, \\
\delta \Omega _{4} &=&\{\mathbf{x}|x_{2}=0,0<x_{1}<R,0<x_{3}<R,\}, \\
\delta \Omega _{5} &=&\{\mathbf{x}|x_{3}=R,0<x_{1}<R,0<x_{2}<R,\}, \\
\delta \Omega _{6} &=&\{\mathbf{x}|x_{3}=0,0<x_{1}<R,0<x_{2}<R,\}.
\end{eqnarray*}
The values of the normal derivatives
$\frac{\partial }{\partial \mathbf{n}}u_{\mathbf{m}}(\mathbf{x})$
of the eigenfunctions on the boundary are equal to certain products
of sine functions:
\begin{equation}
\frac{\partial }{\partial \mathbf{n}}u_{\mathbf{m}}(\mathbf{x})=\left\{
\begin{array}{cccc}
& \frac{8\pi m_{1}}{R^{4}}\sin \frac{\pi m_{2}x_{2}}{R}\sin \frac{\pi
m_{3}x_{3}}{R} & , & \mathbf{x\in }\delta \Omega _{1} \\
(-1)^{m_{1}} & \frac{8\pi m_{1}}{R^{4}}\sin \frac{\pi m_{2}x_{2}}{R}\sin
\frac{\pi m_{3}x_{3}}{R} & , & \mathbf{x\in }\delta \Omega _{2} \\
& \frac{8\pi m_{2}}{R^{4}}\sin \frac{\pi m_{1}x_{1}}{R}\sin \frac{\pi
m_{3}x_{3}}{R} & , & \mathbf{x\in }\delta \Omega _{3} \\
(-1)^{m_{2}} & \frac{8\pi m_{2}}{R^{4}}\sin \frac{\pi m_{1}x_{1}}{R}\sin
\frac{\pi m_{3}x_{3}}{R} & , & \mathbf{x\in }\delta \Omega _{4} \\
& \frac{8\pi m_{3}}{R^{4}}\sin \frac{\pi m_{1}x_{1}}{R}\sin \frac{\pi
m_{2}x_{2}}{R} & , & \mathbf{x\in }\delta \Omega _{5} \\
(-1)^{m_{3}} & \frac{8\pi m_{3}}{R^{4}}\sin \frac{\pi m_{1}x_{1}}{R}\sin
\frac{\pi m_{2}x_{2}}{R} & . & \mathbf{x\in }\delta \Omega _{6}
\end{array}
\right.   \label{bunchofsines}
\end{equation}
As in section 1, in order to reconstruct $f(\mathbf{x})$ we recover Fourier
coefficients $\alpha _{\mathbf{m}}$:
\begin{eqnarray}
\alpha _{\mathbf{m}} &=&\int_{\partial \Omega }I(\mathbf{z},
\lambda _{\mathbf{m}})\frac{\partial }{\partial \mathbf{n}}
u_{\mathbf{m}}(\mathbf{z})ds(\mathbf{z)}  \notag \\
&=&\sum_{j=1}^{6}\int_{\partial \Omega _{j}}I(\mathbf{z},
\lambda _{\mathbf{m}})\frac{\partial }
{\partial \mathbf{n}}u_{\mathbf{m}}(\mathbf{z})ds(\mathbf{z),}
\label{squarecoef}
\end{eqnarray}
where
\begin{equation}
I(\mathbf{z,\lambda })=\int\limits_{0}^{\sqrt{3}R}g(\mathbf{z},r)\Phi
_{\lambda }(r)dr.  \label{wannabeFourier}
\end{equation}
If we choose the Green's function $\Phi _{\lambda }(t)$ in the form
\begin{equation*}
\Phi _{\lambda }(t)=\frac{\cos \lambda t}{4\pi t},
\end{equation*}
equation~(\ref{wannabeFourier}) can be re-written in the form of the Cosine
Fourier transform as follows
\begin{equation}
I(\mathbf{z,\lambda })=\frac{1}{4\pi }\int\limits_{\mathbb{0}}^{\sqrt{3}R}
\left[\frac{g(\mathbf{z},r)}{r}\right] \cos \lambda rdr.  \label{cosinefour}
\end{equation}
As before, when coefficients $\alpha _{\mathbf{m}}$ have been found function
$f(\mathbf{x})$ is obtained by summing the Fourier series
\begin{equation}
f(\mathbf{x})=\sum_{\mathbf{m}\in \mathbb{N}^{3}}\alpha _{\mathbf{m}}
u_{\mathbf{m}}(\mathbf{x}).  \label{newfourierser}
\end{equation}

The above formulae are just a particular case of the inversion technique
presented in the previous section. They yield theoretically exact reconstruction
if the effects of discretization are neglected. However, in the practical
computation only limited range of frequencies
$0\leq \mathbf{\lambda \leq \lambda }^{\mathit{Nyquist}}$ can be recovered
from finitely sampled (in $r$) projections $g(\mathbf{z},r)$ using
equation~(\ref{cosinefour}). Therefore series~(\ref{newfourierser}) has to be
truncated. The Gibbs phenomenon resulting from such a truncation can be
reduced by application of a filter $\eta (\lambda _{\mathbf{m}})$, so that
instead of the previous equation the following formula will be used to
reconstruct an approximation to $f(\mathbf{x}):$
\begin{equation}
f(\mathbf{x})\approx \sum_{\mathbf{m}\in \mathbb{N}^{3},|
\lambda _{\mathbf{m}}|\leq \mathbf{\lambda }^{\mathit{Nyquist}}}\alpha _{\mathbf{m}}\eta
(\lambda _{\mathbf{m}})u_{\mathbf{m}}(\mathbf{x}).  \label{filterfourier}
\end{equation}

The whole reconstruction procedure can be accelerated by utilizing the Fast
Cosine Fourier transform to compute~(\ref{cosinefour}), the 3-D Fast Sine
Fourier transform to sum~series (\ref{newfourierser}), and the 2-D Fast Sine
transform to evaluate the six integrals in equation~(\ref{squarecoef}).
However, there is one obstacle for implementing this plan. The integrals
in~(\ref{squarecoef}) need to be computed for different values of
$\lambda _{\mathbf{m}},$ and there are too many of these values to make
the algorithm efficient. The work-around for this problem is to evaluate
these integrals for a set of uniformly distributed values
$\lambda _{l}=l\Delta \lambda ,$ $l=0,1,2,...$ and then to find the
needed values for each of $\lambda _{\mathbf{m}}$ by interpolation.
Such an interpolation in the spectral
parameter $\lambda $ requires careful selection of discretization steps and
interpolation techniques. The details of our implementation are presented
below.

Suppose function $f(\mathbf{x})$ is to be reconstructed on $n\times n\times
n $ Cartesian grid, and the detectors are located at the nodes of 2-D
$n\times n$ Cartesian grids defined on the faces of the cube
(values at the edges of the cube will not be needed). We will assume that the
discretization step of measurements (in $r)$ is approximately the same as
the step of the Cartesian grids. Then the number of samples $n_{1}$ in one
projection $g(\mathbf{z},r)$ approximately equals $\sqrt{3}n.$ Depending on
a type of the Fast Cosine Fourier transform algorithm used to compute~(\ref
{cosinefour}), the projections will have to be padded with zeros to make the
total number of samples $n_{2}$ either a power of 2, or a product of small
prime numbers. Thus$,$ $n_{1}\leq n_{2}<2n_{1}.$ The Nyquist frequency
$\mathbf{\lambda }^{\mathit{Nyquist}}$ corresponding to such discretization
equals $\pi (n_{1}-1)/D,$ where diameter
$D=\mathrm{diam(}\Omega )=\sqrt{3}R$.
With these parameters in mind we summarize the five steps of the fast
reconstruction algorithm.

\textbf{Step 1}. The first step is to compute a discrete version of (\ref
{wannabeFourier}) using the Fast (discrete) Cosine Fourier transform of
length $n_{2}.$ This will produce values $I(\mathbf{z,\lambda }_{l})$ for
the frequencies
$\lambda _{l}=l\mathbf{\lambda }^{\mathit{Nyquist}}/(n_{2}-1),$
$l=0,1,..,n_{2}-1.$
Importantly, as long as $n_{1}\leq n_{2},$
the step of discretization of $I(\mathbf{z,\lambda }_{l})$ in $\lambda $ is
small enough to make it possible to approximately recover values of
$I(\mathbf{z,\lambda })$ (or a linear function of $I(\mathbf{z,\lambda })$)
for $\lambda \neq \lambda _{l}$ by interpolation.

\textbf{Step 2}. For each value of $\lambda _{l}$ compute integrals in the
form
\begin{equation}
\int_{\partial \Omega _{j}}I(\mathbf{z},\lambda _{l})\sin \frac{\pi
m_{i}x_{i}}{R}\sin \frac{\pi m_{k}x_{k}}{R}ds(\mathbf{z)},~j=1,...6,
\label{fourierintegrals}
\end{equation}
for integer all integer $m_{i},$ $m_{k}$ by means of the 2-D Fast Sine
Fourier transform.

\textbf{Step 3}. Compute approximate values of integrals
\begin{equation*}
\int_{\partial \Omega _{j}}I(\mathbf{z},\lambda _{\mathbf{m}})\sin \frac{\pi
m_{i}x_{i}}{R}\sin \frac{\pi m_{k}x_{k}}{R}ds(\mathbf{z)},~j=1,...6,
\end{equation*}
by interpolating values obtained on step 2
(equation \ref{fourierintegrals}).
Some care should be taken to guarantee accuracy of computations on this
step. It is well known that a low order interpolation in spectral parameter
can lead to a suboptimal reconstruction, due to highly oscillatory nature of
the Fourier transformant. An example and analysis of this phenomenon can be
found in \cite{Natt}. However, the Fourier transform of a finitely supported
function is an analytic function of the spectral parameter, even if the
function itself is known imprecisely. Therefore, higher order polynomial
interpolation is applicable and does produce good results in this case. In
our numerical experiments we observed that, indeed, the linear interpolation
in $\lambda $ yields rather inaccurate reconstruction. The increase in the
order of the polynomial interpolation significantly improves the image; if
the 6-th order Lagrange interpolation is utilized, the interpolation error
is dominated by the discretization errors and further increase in the
accuracy of interpolation is not needed.

\textbf{Step 4}. Use values computed on step 3 to calculate
$\lambda _{\mathbf{m}}$
by combining equations (\ref{bunchofsines}) and (\ref{squarecoef}).

\textbf{Step 5}. Using the 3-D Fast Fourier Sine transform to implement~(\ref
{filterfourier}), compute values of $f(\mathbf{x})$ at the nodes of the 3-D
Cartesian grid.

A simple computation shows that the number of floating point operations
implemented on each steps of the algorithm is $\mathcal{O}(n^{3}\log n)$ for
steps 1,2, and 5, and $\mathcal{O}(n^{3})$ on steps 3 and 4, resulting in a
total $\mathcal{O}(n^{3}\log n)$ operation count for this technique. This
has to be compared with $\mathcal{O}(n^{5})$ operation count required by a
backprojection step of a method resulting from a straightforward
discretization of any of the explicit inversion
formulae\cite{Finch,XuWang,kunya}.
(The latter estimate assumes that the reconstruction
is done on $n\times n\times n$ Cartesian grid from $n^{2}$ detector
positions.)

In what follows we present a numerical example illustrating the work of our
algorithm. A function is reconstructed within the cube $[0,1]\times [
0,1]\times [0,1].$ The dimension of the grids were defined by the
values of parameters $n=129,$ $n_{1}=223,$ and $n_{2}=256.$ As a filter we
used the cosine window function
\begin{equation*}
\eta (\lambda )=\left\{
\begin{array}{cc}
\cos \frac{\pi \lambda }{2\lambda ^{\mathit{Nyquist}}}, & \lambda \leq
\lambda ^{\mathit{Nyquist}} \\
0, & \lambda >\lambda ^{\mathit{Nyquist}}
\end{array}
\right. .
\end{equation*}
We utilized the same phantom as in~\cite{kunya} to facilitate the comparison
of the present results with the images reconstructed in the former work by
application of discretized explicit inversion formulae. The phantom consists
of eight characteristic functions of the balls with radii ranging from 0.06
to 0.13, whose centers lie in the plane $x_{3}=0.$ The cross section of the
phantom by the latter plane is shown in Figure~1(a). Figure~1(b) shows the
central cross section ($x_{3}=0$) of the reconstruction from the exact data.
\begin{figure}[th]
\begin{center}
\begin{tabular}{ccc}
\includegraphics[width=2.0in,height=2.0in]{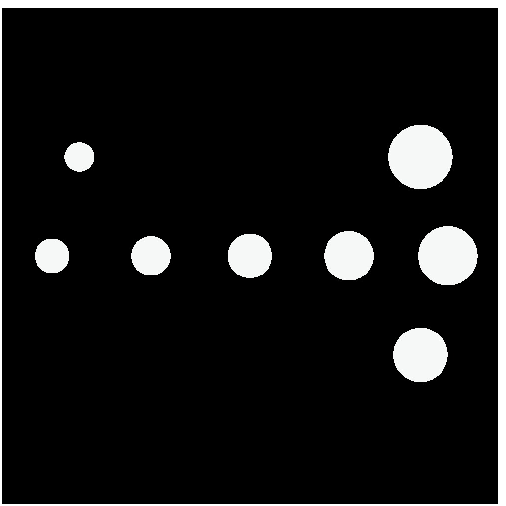} &
\includegraphics[width=2.0in,height=2.0in]{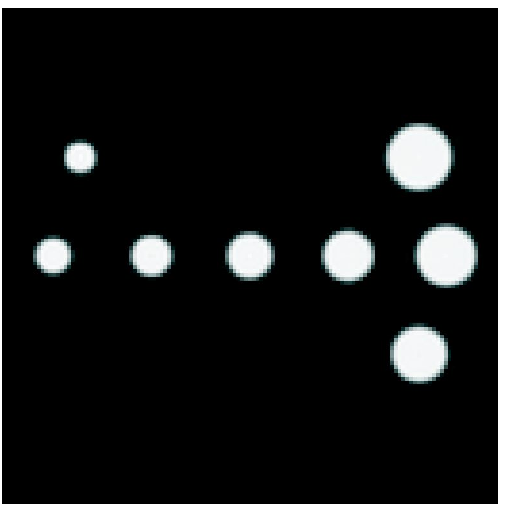} &
\includegraphics[width=2.0in,height=2.0in]{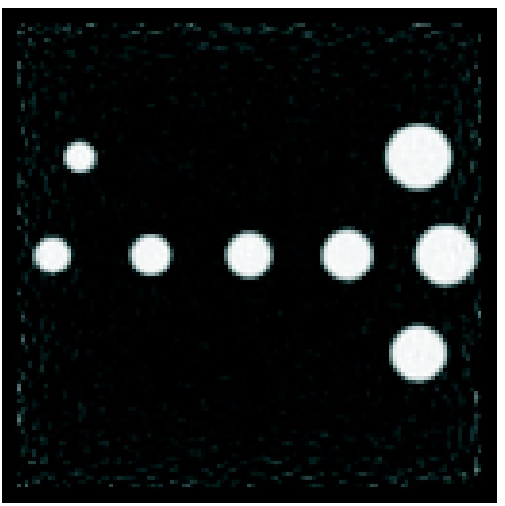} \\
&  &  \\
(a) & (b) & (c)
\end{tabular}
\label{fig1}
\end{center}
\caption{Numerical example: (a) the phantom, (b) reconstruction from the
exact data, and (c) reconstruction from noisy data}
\end{figure}

In order to evaluate the sensitivity of the algorithm to noise in data, the
imprecise measurements were modeled by adding to the projections normally
distributed noise with the intensity 15\% of the signal (in $L_{2}$-norm).
In this experiment the values of the reconstructed function were set to zero
outside of the cube $[0.05,0.95]\times [0.05,0.95]\times [
0.05,0.95]$, since the reconstruction from noisy data is unstable at the
locations close to the detectors, due to the singular nature of the Green's
function. Such sensitivity is natural, and is an issue for other
reconstruction techniques as well. For example, the slices of 3-D images
obtained in~\cite{kunya} were computed within the unit ball from the
detectors located on a sphere of radius $1.1;$ the reconstruction in the
close vicinity of the detectors was also avoided. The gray scale
representation of the reconstructed function is shown in Figure~1(c); in
Figure~2 we demonstrate the surface plot of the same function.
\begin{figure}[h]
\begin{center}
\includegraphics[width=4.78in,height=2.2in]{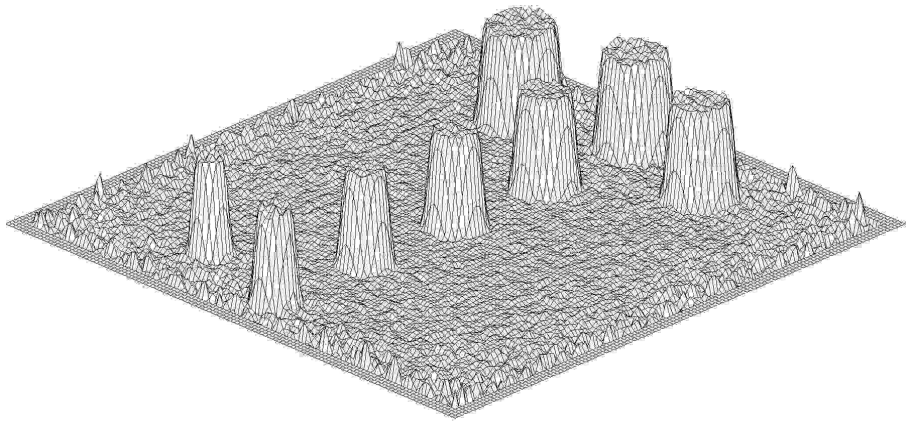} \label{aksnoise}
\end{center}
\caption{Reconstruction from noisy data; surface plot}
\end{figure}

We believe that the quality of the reconstructed images in all of the above
experiments is as good as of those resulting from explicit inversion
formulae (see~\cite{kunya}). Similarly to those formulae the present
technique demonstrates high stability of the algorithm to the perturbations
of the data. On the other hand, the computation time for the present
algorithm in the above experiment varied from 7 to 8 seconds on an
AMD workstation with a 2GHz processor; function values at about 2 million points
were reconstructed from about 97 thousand projections.
For comparison, the reconstruction reported in~\cite{kunya} from 33000 projections
at about a million grid points by the fastest of our implementations of the
explicit inversion formula took about 48 minutes. If the latter algorithm is
used to process the same number of the projections and grid points as
in the present example (97 thousands and 2 million respectively), the
computation time increases to about 7 hours. It is fair to say that
the fast algorithm is thousands time faster than the straightforward
discretization of any of the backprojection-type formulae.

\section{Reconstruction in the presence of exterior sources}

The series solution described above has an interesting property not
possessed (to the best of our knowledge) by any other currently known
explicit reconstruction technique. Let us consider a slightly more general
problem. Suppose that region $\Omega $ is a proper subset of a larger region
$\Omega _{1}$ ($\Omega \subset \Omega _{1}$) and that a $L^{2}$ function $F$
is defined on $\Omega _{1}.$ We will denote the restriction of $F$ on
$\Omega $ by $f,$ i.e.
\begin{equation*}
f(\mathbf{x})=\left\{
\begin{array}{ccc}
F(\mathbf{x}) & , & \mathbf{x}\in \Omega  \\
0 & , & \mathbf{x}\in \mathbb{R}^{n}\backslash \Omega
\end{array}
\right. .
\end{equation*}
We would like to reconstruct $f(\mathbf{x})$ from the integrals
$g(\mathbf{z},r)$ of $F$ over spheres with the centers
on $\partial \Omega $:
\begin{equation*}
g(\mathbf{z},r)=\int\limits_{\mathbb{S}^{n-1}}
F(\mathbf{z}+r\mathbf{\hat{s}})r^{n-1}
d\mathbf{\hat{s}},\qquad \mathbf{z}\in \partial \Omega .
\end{equation*}
Unlike in the previously considered problem, now the centers of the
integration spheres are lying on a surface contained within the support
$\Omega _{1}$ of the function $F.$ While we are still trying to reconstruct
the restriction $f$ of $F$ to $\Omega $, the integrals we know are those of
$F$ and not of $f$.
\begin{figure}[th]
\begin{center}
\begin{tabular}{cc}
\includegraphics[width=2.0in,height=2.0in]{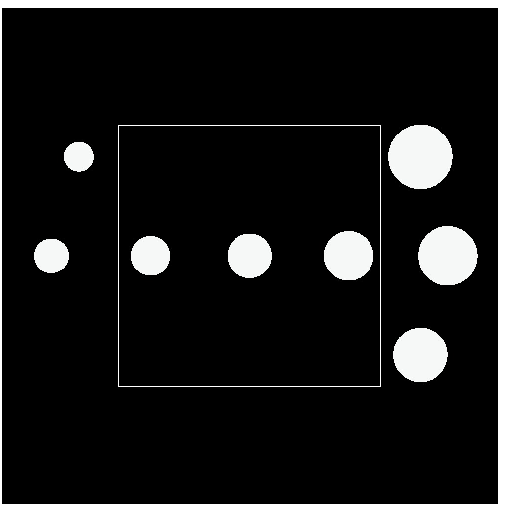} &
\includegraphics[width=2.0in,height=2.0in]{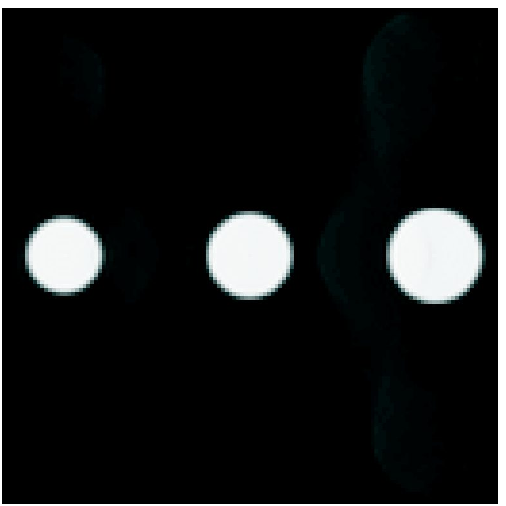} \\
&  \\
(a) & (b)
\end{tabular}
\label{mappart}
\end{center}
\caption{Reconstruction in the presence of exterior sources (a) phantom; the
white line shows location of the detectors (b) reconstructed image.}
\end{figure}

It turns out that the solution to this problem is still given by
formulae~(\ref{fourierser}) and~(\ref{serkoef1a})
(or equivalently by~(\ref{serbackpr}), (\ref{serfiltr}),
and~(\ref{serkoef1a})).
Indeed, if we extend functions
$u_{m}(\mathbf{x})$ by $0$ to $\mathbb{R}^{n}$, formula~(\ref{helmdiscr})
holds for all $\mathbf{x}\in \mathbb{R}^{n}$, and~(\ref{fourierser}) remains
unchanged. In formula~(\ref{serkoef1}) $f$ can be replaced by $F$ as
follows:
\begin{align*}
\alpha _{m}& =\int_{\Omega }u_{m}(\mathbf{x})f(\mathbf{x})d\mathbf{x=}
\int_{\Omega }u_{m}(\mathbf{x})F(\mathbf{x})d\mathbf{x} \\
& =\int_{\partial \Omega }\left( \int_{\Omega }\Phi _{\lambda _{m}}
(|\mathbf{x}-\mathbf{z}|)F(\mathbf{x})d\mathbf{x}\right)
\frac{\partial }{\partial
\mathbf{n}}u_{m}(\mathbf{z})ds(\mathbf{z),}
\end{align*}
and the inner integral can be computed from projections as before:
\begin{equation*}
\int_{\Omega }\Phi _{\lambda _{m}}(|\mathbf{x}-\mathbf{z}|)F(\mathbf{x})
d\mathbf{x}=\int\limits_{\mathbb{R}^{+}}g(\mathbf{z},r)\Phi _{\lambda
_{m}}(r)dr.
\end{equation*}
By combining the two above equations we again arrive at the
formula~(\ref{serkoef1a}).

This interesting property can be illustrated by a numerical example. We
consider the same phantom as in the previous section. This time, however,
the detectors are located on the surface of the cube $\Omega
=[0.235,0.765]\times [0.235,0.765]\times [0.235,0.765].$
Location of the detectors is shown in Figure~3(a) by a white line. The
integrals were computed over full spheres, and the reconstruction was
conducted, as before on the grid of size $129\times 129\times 129$ within
$\Omega$. The result is presented in Figure~3(b), and, as a surface plot, in
Figure~4. On the latter figure one can notice shallow troughs resulting from
imperfect resolution of sharp edges of exterior sources by a finite number
of detectors. The depths of these troughs, however, does not exceed $6\%$ of
the maximum of the original function.
\begin{figure}[h]
\begin{center}
\includegraphics[width=5.0in,height=2.2in]{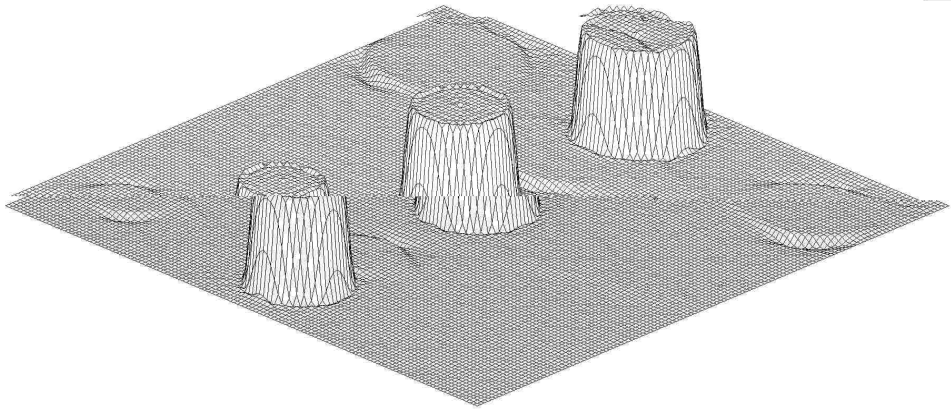} \label{akspart}
\end{center}
\caption{Reconstruction in the presence of exterior sources; surface plot.}
\end{figure}

\section{Acknowledgments}

The author would like to thank P. Kuchment for fruitful discussions and
numerous helpful comments.

This work was partially supported by the NSF/DMS grant NSF-0312292 and by
the DOE grant DE-FG02-03ER25577.


\begin{thebibliography}{99}
\bibitem{trinity}  M. Agranovsky, P. Kuchment, and E. T. Quinto, Range
descriptions for the spherical mean Radon transform, preprint 2006,
\emph{arXiv: math}. AP/0606314.

\bibitem{AQ}  M. L. Agranovsky and E. T. Quinto, Injectivity sets for the
Radon transform over circles and complete systems of radial functions,
\emph{Journal Of Functional Analysis}, \textbf{139}, no. 2 (1996) pp. 383-414.

\bibitem{AK1}  G. Ambartsoumian and P. Kuchment, On the injectivity of the
circular Radon transform arising in thermoacoustic tomography,
\emph{Inverse Problems} \textbf{21} (2005), pp. 473-485.

\bibitem{AK2}  G. Ambartsoumian and P. Kuchment, A range description for the
planar circular Radon transform, \emph{SIAM J. Math. Anal.} \textbf{38}, no.
2 (2006) pp. 681-692.



\bibitem{berard1}  P. B\'{e}rard, Spectres et Groupes Cristallographiques,
\emph{\ C. R. Acad. Sci. Paris A-B}\textbf{\ 288,} no. 23 (1979) pp.
A1059-A1060.

\bibitem{berard2}  P. B\'{e}rard and G. Besson, Spectres et Groupes
Cristallographiques II: Domaines Sph\'{e}riques, \emph{Ann. Institut Fourier},
\textbf{30,} no. 3 (1980) pp. 237-248.

\bibitem{Faw} J. A. Fawcett, Inversion of $N$-dimensional spherical averages,
\emph{SIAM J. Appl. Math.}  \textbf{45}, no. 2, (1985) pp. 336-341.

\bibitem{Finch}  D. Finch, Rakesh, and S. Patch, Determining a function from
its mean values over a family of spheres, \emph{SIAM J. Math. Anal.}
\textbf{35} no. 5 (2004), 1213--1240.

\bibitem{Finch1}  D. Finch and Rakesh, The range of the spherical mean value
operator for functions supported in a ball, \emph{Inverse Problems}
\textbf{22} (2006), pp. 923-938.

\bibitem{Finch2}  D. Finch, M. Haltmeier, and Rakesh, Inversion of spherical
means and the wave equation in even dimensions, preprint 2006.

\bibitem{kuchrev}  P. Kuchment, Generalized Transforms of Radon Type and
Their Applications, in G. Olafsson and E. T. Quinto (Editors), \emph{The
Radon Transform, Inverse Problems, and Tomography,} Proc. Symp. Appl. Math.
v. 63, AMS, Providence, RI 2006, pp.67 - 91.

\bibitem{kunya}  L. A. Kunyansky, Explicit inversion formulae for the
spherical mean Radon transform, to appear in \emph{Inverse Problems}.

\bibitem{kruger1}  R. A. Kruger, P. Liu, Y. R. Fang, and C. R. Appledorn,
Photoacoustic ultrasound (PAUS) reconstruction tomography, \emph{Med. Phys}.
\textbf{22,} (1995), pp. 1605-1609.

\bibitem{kruger}  R. A. Kruger, D. R. Reinecke, and G. A.Kruger, GA
Thermoacoustic computed tomography - technical considerations, \emph{Med.
Phys.} \textbf{26,} no. 9, (1999) pp. 1832-1837.

\bibitem{Natt}  F. Natterer and F. Wubbeling, Mathematical Methods in Image
Reconstruction, \emph{Monographs on Mathematical Modeling and Computation},
\textbf{5}, SIAM, Philadelphia, PA 2001.

\bibitem{Norton2D}  S. J. Norton, Reconstruction of a two-dimensional
reflecting medium over a circular domain: Exact solution, \emph{J. Acoust.
Soc. Amer.} \textbf{67} (1980), pp. 1266--1273.

\bibitem{Nort3D}  S. J. Norton and M. Linzer, Ultrasonic Reflectivity
Imaging in Three Dimensions - Exact Inverse Scattering Solutions for Plane,
Cylindrical, and Spherical Apertures, \emph{IEEE Trans. Biomed. Eng.}
\textbf{28}, no. 2 (1981) pp. 202-220.

\bibitem{XuWang0}  M. Xu and L.-H. V. Wang, Time-domain reconstruction for
thermoacoustic tomography in a spherical geometry, \emph{IEEE Trans. Med.
Imag.} \textbf{21,} (2002), pp. 814-822.

\bibitem{XuWang}  M. Xu, L. V. Wang, Universal back-projection algorithm for
photoacoustic computed tomography, \emph{Phys Review E} \textbf{71} (2005),
p. 016706.
\end{thebibliography}
\end{document}